\documentclass{amsart}
\usepackage{amssymb,amscd}
\usepackage{graphicx}

\usepackage{latexsym}

\usepackage{maple2e}

\DefineParaStyle{Maple Normal}

\usepackage{amsthm}
\newtheorem{theorem}{Theorem}
\newtheorem{lemma}{Lemma}

\newcommand{\bc}{\mathbb{C}}
\newcommand{\br}{\mathbb{R}}
\newcommand{\bp}{\mathbb{P}}
\newcommand{\bq}{\mathbb{Q}}
\newcommand{\modm}{\mathcal{M}}
\newcommand{\modt}{\mathcal{T}}
\newcommand{\ay}{\mathcal{A}}
\newcommand{\bb}{\mathcal{B}}

\begin{document}

\title{Weil-Petersson volumes and cone surfaces}
\author{Norman Do and Paul Norbury}
\address{Department of Mathematics and Statistics\\
University of Melbourne\\Australia 3010}
\email{N.Do@ms.unimelb.edu.au,\ pnorbury@ms.unimelb.edu.au}
\thanks{The second author would like to thank MSRI for support during completion of this work}

\date{\today}

\begin{abstract}

\noindent The moduli spaces of hyperbolic surfaces of genus $g$ with $n$ geodesic boundary components are naturally symplectic manifolds.  Mirzakhani proved that their volumes are polynomials in the lengths of the boundaries by computing the volumes recursively.  In this paper we give new recursion relations between the volume polynomials.

\end{abstract}

\maketitle

\section{Introduction}

For ${\bf L}=(L_1,L_2,...,L_n)$, a sequence of non-negative numbers, let  $\modm_{g,n}({\bf L})$ be the moduli space of connected oriented genus $g$ hyperbolic surfaces with $n$ labeled boundary components of lengths $L_1,...,L_n$.  A cusp at a point naturally corresponds to a zero length boundary component.  When ${\bf L}=0$, that is there are $n$ cusps, the moduli space  $\modm_{g,n}({\bf 0})$ is naturally identified with the moduli space of conformal structures on a genus $g$ oriented surface with $n$ labeled points, also known as the moduli space of curves with $n$ labeled points.  The identification uses the fact that in any conformal class of metrics there is a unique complete hyperbolic metric, and for every conformal automorphism there is a corresponding isometry.

On the moduli space $\modm_{g,n}({\bf L})$ lives a natural symplectic form $\omega$, defined precisely in Section~\ref{sec:vol}.  The volume of the moduli space is
\[
V_{g,n}({\bf L})=\int_{\modm_{g,n}({\bf L})}\frac{\omega^{3g-3+n}}{(3g-3+n)!},\quad (g,n)\neq(1,1).
\]
When $(g,n)=(1,1)$ we instead take half of the integral of $\omega$, an orbifold volume,
\[ V_{1,1}(L_1)=\frac{1}{2}\cdot\int_{\modm_{1,1}(L_1)}\omega=\frac{1}{48}(L_1^2+4\pi^2),\]
which fits well with recursion relations between volumes, and relations with intersection numbers on the moduli space.  Mirzakhani uses the true volume of $\modm_{1,1}(L_1)$ in \cite{MirSim,MirWei} and includes an extra factor of a half in her formulae.
\begin{theorem}[Mirzakhani \cite{MirSim}]  \label{th:mirz}
$V_{g,n}({\bf L})$ is a polynomial in ${\bf L}=(L_1,...,L_n)$.  The coefficient of $L^{\alpha}=L^{\alpha_1}..L^{\alpha_n}$ lies in $\pi^{6g-6+2n-|\alpha|}\bq$, $|\alpha|=\alpha_1+...+\alpha_n$.
\end{theorem}
Mirzakhani proved this using a recursion relation between volumes of moduli spaces:  
\begin{equation} \label{eq:rec}
\frac{\partial}{\partial L_1} (L_1V_{g,n}({\bf L}))=\ay_{g,n}({\bf L})+\bb_{g,n}({\bf L})
\end{equation}
where $\ay_{g,n}({\bf L})$ consists of integral transforms of $V'_{g-1,n+1}$ and $\bb_{g,n}({\bf L})$ consists of integral transforms of $V_{g,n-1}$.   The $'$ in $V'_{g-1,n+1}$ allows disconnected surfaces which is an efficient way to encode the volumes of pairs of moduli spaces, or equivalently volumes of products $\modm_{g_1,n_1}\times \modm_{g_2,n_2}$ for $g_1+g_2=g$ and $n_1+n_2=n+2$.  We have omitted the ${\bf L}$ dependence in $V'_{g-1,n+1}$ and $V_{g,n-1}$ because it requires further explanation.  See Section~\ref{sec:coord} for precise definitions of $\ay_{g,n}({\bf L})$ and $\bb_{g,n}({\bf L})$. 

The main idea of this paper is to use intermediary moduli spaces to give new recursion relations between volumes of moduli spaces.  The intermediary moduli spaces consist of hyperbolic surfaces with a cone point of a specified angle.  Hyperbolic geometry is an ideal setting for studying cone points, although a cone point does make sense more generally in terms of a conformal structure on a Riemann surface.  A cone angle of 0 corresponds to a cusp marked point and as the cone angle goes from 0 to $2\pi$ this corresponds, in a sense, to removing the marked point.  This leads to interesting relations between the moduli spaces.  These intermediary moduli spaces are reminiscent of the moduli spaces of anti self dual connections with cone singularities around an embedded surface in a four-manifold, used by Kronheimer and Mrowka \cite{KMroEmb} to get relationships between intersection numbers on instanton moduli spaces.   

In \cite{TWZGen} it is shown that one can interpret a point with cone angle in terms of an imaginary length boundary component.  Explicitly, a cone angle $\phi$ appears by substituting the length $i\phi$ in the volume polynomial.  Mirzakhani's results, Theorem~\ref{th:mirz} and (\ref{eq:rec}) use a generalised McShane formula \cite{McSRem} on hyperbolic surfaces, which was adapted in \cite{TWZGen} to allow a cone angle $\phi$ that ends up appearing as a length $i\phi$ in such a formula, and hence in the volume polynomial.  We do not describe the generalised McShane formula in this paper, although in Section~\ref{sec:coord} we give the underlying idea in terms of coordinates on the hyperbolic surface.
\begin{theorem}  \label{th:string}
For ${\bf L}=(L_1,...,L_n)$
 \begin{equation}  \label{eq:string}
 V_{g,n+1}({\bf L},2\pi i)=\sum_{k=1}^n\int_0^{L_k}L_kV_{g,n}({\bf L})dL_k
 \end{equation}
and
\begin{equation}  \label{eq:dilaton}
\frac{\partial V_{g,n+1}}{\partial L_{n+1}}({\bf L},2\pi i)=2\pi i(2g-2+n)V_{g,n}({\bf L}).
\end{equation}
\end{theorem}
We think of the theorem as describing the limit of the volume and its derivative when a cone angle tends to $2\pi$, and hence is removable, although the statement of the theorem is independent of this interpretation.  The recursion relations (\ref{eq:string}) and (\ref{eq:dilaton}) give information about the volume of $\modm_{g,n+1}(L_1,...,L_{n+1})$ from the volume of the lower dimensional moduli space $\modm_{g,n}(L_1,...,L_n)$.  This contrasts with Mirzakhani's relation (\ref{eq:rec}) which uses many lower dimensional moduli spaces as described above.  We discuss this difference more in Section~\ref{sec:rec} and in particular give a simpler algorithm to determine $V_{0,n}({\bf L})$.

There are three potential approaches to the proof of Theorem~\ref{th:string}, one of which we present here, in Section~\ref{sec:proof}, and the others we cannot yet do.  Firstly, the theorem must follow from Mirzakhani's recursion relation (\ref{eq:rec}) since that relation uniquely determines the polynomials.  Secondly, the theorem should follow from an analysis of the cone angle moduli spaces.  Thirdly, the statement of the theorem is equivalent to relations between the coefficients of the volume polynomials which are intersection numbers of $\psi$ classes and $\kappa$ classes (see Section~\ref{sec:proof} for definitions), so relations between the latter can be used to deduce the theorem.  It is this third approach that we present here.

\section{Volume of the moduli space}   \label{sec:vol}
\subsection{Fenchel-Nielsen coordinates.}
Given an oriented hyperbolic surface $\Sigma$ with $n$ geodesic boundary components, cut along a maximal set of disjoint embedded simple closed geodesics.  The resulting pieces are genus zero hyperbolic surfaces each with 3 geodesic boundary components, i.e. hyperbolic pairs of pants. 
Each pair of pants contributes Euler characteristic -1, so there are $2g-2+n=-\chi(\Sigma)$ pairs of pants
in the decomposition, and hence $3g-3+n$ closed geodesics (not counting the boundary geodesics.) We call such a decomposition a {\em pair of pants decomposition} of the surface $\Sigma$.   A pair of pants decomposition gives a coordinate system on the moduli space of hyperbolic surfaces - the lengths of the disjoint embedded simple closed geodesics and angles which we will now define.

Each boundary component of an oriented hyperbolic pair of pants has a pair of distinguished points - the closest points to the other two boundary components - which divide the boundary circle into two equal length pieces.  A closed geodesic in a pair of pants decomposition contains two pairs of distinguished points coming from each side of the geodesic.  The (oriented) angle between the pairs is a number $\theta\in[0,l]$, where $l$ is the length of the geodesic, well-defined up to $\theta\mapsto\theta+l/2$.

So we have $(l_i,\theta_i)$ for $i=1,...,3g-3+n$, the Fenchel-Nielsen coordinates.  Although $\theta_i\in[0,l_i)$ it is natural instead to lift it to $\theta_i\in\br$.  The coordinates take values in 
\[\modt_{g,n}({\bf L})=(\br^+\times\br)^{3g-3+n}\]
which is known as {\em Teichm\"uller space}.  This statement also uses the fact that any triple of non-negative numbers give boundary lengths of a unique hyperbolic pair of pants and two oriented hyperbolic pairs of pants can be glued along boundaries of common length with any angle between distinguished points.

The coordinates are canonical coordinates for a symplectic form
\begin{equation}  \label{eq:symp}
\omega=\sum_i dl_i\wedge d\theta_i
\end{equation}
which is a rather trivial statement.  What is quite deep is the fact that the symplectic form is invariant under the action of the mapping class group $Mod_{g,n}$, of isotopy classes of orientation preserving diffeomorphisms of the surface that preserve boundary components, and hence descends to  a quotient of Teichm\"uller space 
\[ \modm_{g,n}({\bf L})=\modt_{g,n}({\bf L})/Mod_{g,n}\]
known as the moduli space.  The action can be understood by noticing that a given hyperbolic surface has infinitely many pair of pants decompositions that are geometrically different (but topologically the same and thus related by an element of $Mod_{g,n}$.)  Each different decomposition gives different lengths and angles to the same hyperbolic surface, and hence different coordinates, whereas the symplectic form (\ref{eq:symp}) is preserved.  This gives rise to the volume of the moduli space $V_{g,n}({\bf L})$ given by the integral of $\omega^{3g-3+n}/(3g-3+n)!$ over the moduli space, or equivalently over a fundamental domain for $Mod_{g,n}$ in $\modt_{g,n}({\bf L})$.

\subsection{Coordinates on a hyperbolic surface}  \label{sec:coord}
It is useful to view a hyperbolic surface from one boundary component chosen from the $n$ boundary components.  To motivate this idea, let us consider isometries of the surface that leave the boundary components invariant, although they may be rotated.   A non-trivial isometry detects when the mapping class group does not act freely on $\modt_{g,n}({\bf L})$.  An important fact is that {\em an isometry that fixes any boundary component is necessarily the identity.}  To see this, note that any point on the surface lies on a geodesic that meets the $i$th boundary component, say, at right angles - simply take the shortest path to that boundary component.  (The path lies on the interior of the surface since if it was to meet another boundary component, it would meet at a positive angle - two geodesics cannot be tangent - and thus there would be a shorter path by taking a shortcut to the boundary near the angle.)  An isometry that fixes the $i$th boundary component will preserve the orthogonal geodesic, so in particular it fixes all points along that geodesic.  Applying this argument to any point on the surface shows that all points of the surface are fixed by the isometry.   Thus any specified boundary component, say $\partial_1$, detects non-trivial isometries.

Another way to view the discussion in the previous paragraph is that the geodesic boundary component  $\partial_1$ gives a coordinate system on the surface - to any point on the surface assign its distance from $\partial_1$  and the point on $\partial_1$  where the shortest geodesic meets.  
More generally, take any geodesic beginning at a given point on the surface and meeting $\partial_1$  perpendicularly, and assign to the point its length and the point it meets $\partial_1$ .  This
makes the coordinate system locally smooth, at the cost of losing uniqueness for the coordinates of a point.

Mirzakhani uses this coordinate system in the following way.  Project points onto the second coordinate, which takes its values in $\partial_1$.  Now suppose that there is another boundary component, $\partial_i$ say.  The projection of $\partial_i$ is an interval $I^0_i\subset \partial_1$.  More precisely, the projection is a collection of infinitely many disjoint intervals $\{ I^j_i\ |\ j=0,...,\infty\}$ since we take any perpendicular geodesic, not just the shortest one, resulting in non-unique coordinates.

The sum of the lengths $f_i=\sum_jl( I^j_i)$ is a well-defined function on the moduli space ${\modm}_{g,n}({\bf L})$.  The length of a single interval $l( I^j_i)$ is well-defined on Teichm\"uller space ${\modt}_{g,n}({\bf L})$, and although it does not descend to the moduli space, $l( I^j_i)$ descends to an intermediate moduli space:
\[\begin{array}{c}
\modt_{g,n}({\bf L})\\\downarrow\\\widehat{\modm}_{g,n}({\bf L})\\\downarrow\\\modm_{g,n}({\bf L})
\end{array}\]
and Mirzakhani shows that this enables one to integrate the function $f_i=\sum_jl( I^i_j)$ over $\modm_{g,n}({\bf L})$ yielding a polynomial, calculable from $V_{g,n-1}$.  The $n-1$ collections of intervals $\{ I^j_i\ |\ j=0,...,\infty\}$, $i=2,...,n$ are disjoint from each other and Mirzakhani similarly shows that the complementary region (up to a measure zero set) gives a well-defined function $f^c$ on the moduli space which can be integrated in terms of lower volumes.    
Since $f^c+\sum f_i=L_1$, the sum of all of the integrals gives
\[\int_{\modm_{g,n}({\bf L})}L_1d{\rm vol}=L_1V_{g,n}({\bf L})\]
the derivative of which can be calculated and leads to Mirzakhani's recursion relation:
\[\frac{\partial}{\partial L_1} (L_1V_{g,n}({\bf L}))=\ay_{g,n}({\bf L})+\bb_{g,n}({\bf L}).\]
For completeness we will define the right hand side although this will not be used further in the paper.
Put  ${\hat{\bf L}}=(L_2,...,L_n)$ and let $(L_2,...,\hat{L}_j,...,L_n)$ mean we remove $L_j$.  Then
\[\ay_{g,n}({\bf L})=\int K_{L_1}(x,y)V'_{g-1,n+1}(x,y, {\hat{\bf L}})dxdy\]
where 
\[ V'_{g-1,n+1}(x,y, {\hat{\bf L}})=V_{g-1,n+1}(x,y, {\hat{\bf L}})+\sum_{g_i,n_i,{\bf L_i}}
V_{g_1,n_1}(x,{\bf L_1})\times V_{g_2,n_2}(y,{\bf L_2})\]
and the sum is over all $g_1+g_2=g$,  $n_1+n_2=n+1$ and ${\bf L_1}\sqcup{\bf L_2}={\hat{\bf L}}$.
And
\[\bb_{g,n}({\bf L})=\sum_{j=2}^n\int K_{L_1,L_j}(x)V_{g,n-1}(x, L_2,...,\hat{L}_j,...,L_n)dx.\]
The kernels are defined by 
\[ K_{L_1}(x,y)=H(x+y,L_1),\  \ \ K_{L_1,L_j}(x)=H(x,L_1+L_j)+H(x,L_1-L_j)\]
for
\[ H(x,y)=\frac{1}{2}\left(\frac{1}{1+e^{\frac{x+y}{2}}}+\frac{1}{1+e^{\frac{x-y}{2}}}\right).\]
The derivation of these kernels comes from a detailed study of a hyperbolic pair of pants - the simplest hyperbolic surface to contain two boundary components.  We refer the reader to \cite{MirSim,MirWei} for full details.

\section{Characteristic classes of surface bundles}  \label{sec:proof}
\subsection{Surface bundles.}
To any oriented topological surface bundle
\[ \begin{array}{ccrll}\Sigma_g&\hookrightarrow&X&\\&&\pi\downarrow&\uparrow s_i\\&&B&\end{array}
\quad\quad i=1,...,n\]
with $n$ sections having disjoint images we can associate characteristic classes in $H^*(B)$, \cite{MorCha}. 
On $X$ there is a complex line bundle $\gamma\to X$ with fibre at $b\in B$ the vertical cotangent bundle $T^*\pi^{-1}(b)$.  A local trivialisation is obtained from a local trivialisation of the fibre bundle $X$.  For each $i=1,...,n$ pull back the line bundle $\gamma$ to $s_i^*\gamma=\gamma_i\to B$.  Define 
\[ \psi_i=c_1(\gamma_i)\in H^2(B).\]
Let $e=c_1(\gamma)\in H^2(X)$.  (We use the terminology $e$ because it is naturally the Euler class of $\gamma$.  We have put a complex structure on $\gamma$ for convenience.)  Define the Mumford-Morita-Miller classes
\[\tilde{\kappa}_m=\pi_!e^{m+1}\in H^{2m}(B)\] 
where $\pi_!:H^k(X)\to H^{k-2}(B)$ is the umkehr map, or Gysin homomorphism, obtained by integrating along the (oriented) fibres.  Alternatively, the umkehr map is obtained from the composition
\[ \pi_!:H^k(X)\stackrel{PD}{\to}  H_{d-k}(X)\stackrel{\pi_*}{\to} H_{d-k}(B)\stackrel{PD}{\to} H^{k-2}(B)\]
where $d=\dim X$ and $PD$ denotes Poincare duality.
The Mumford-Morita-Miller classes ignore the $n$ sections $s_i$.  Use instead the sequence 
\[ \pi_!:H^k_c(X-\cup s_i(B))\stackrel{PD}{\to}  H_{d-k}(X-\cup s_i(B))\stackrel{\pi_*}{\to}  H_{d-k}(B)\stackrel{PD}{\to} H^{k-2}(B)\]
where $H^k_c(X-\cup s_i(B))$ denotes cohomology with compact supports.  Define the kappa classes
\[\kappa_m=\pi_!e_c^{m+1}\in H^{2m}(B)\]
where $e_c=e(\gamma)\in H^2_c(X-\cup_i s_i(B))$ is the Euler class with compact support.  It has the property that on any fibre $\Sigma$
\[ \langle e_c,\Sigma-\cup_i s_i(B)\rangle=-\chi(\Sigma-\cup_i s_i(B))\]
which generalises $\langle e,\Sigma\rangle=-\chi(\Sigma)$.  It is convenient to work with the compact manifold $X$ and in place of $e_c$ use its image in $H^2(X)$
\begin{eqnarray*}
H^2_c(X-\cup_i s_i(B))&\to& H^2(X)\\
e_c&\mapsto& e_n=e+\sum_{i=1}^n PD[s_i(B)].
\end{eqnarray*}
The expression for $e_n$ is deduced from its two properties
\begin{equation}  \label{eq:euler}
\langle e_n,\Sigma\rangle=-\chi(\Sigma-\cup_i s_i(B)),\quad e_n\cdot PD[s_j(B)]=0\ \ j=1...n
\end{equation}
the first because it is defined by restriction, and the second because it lies in the kernel of the map
$H^2(X)\to H^2(\cup_i s_i(B))$.

We will need relations between classes obtained by simply forgetting a section.  Now
\[ e_{n+1}=e_n+PD[s_{n+1}(B)]\]
so from $e_{n+1}\cdot PD[s_{n+1}(B)]=0$ and $c=a+b\Rightarrow c^{m+1}=a^{m+1}+b\sum_{j=0}^mc^ja^{m-j}$
\[ e_{n+1}^{m+1}=e_n^{m+1}+PD[s_{n+1}(B)]\cdot e_n^m\]
thus the forgetful map $\pi_{n+1}$ induces $\pi_{n+1}^*:H^*(B)\to H^*(B)$ satisfying
\begin{equation}  \label{eq:pbkap}
\kappa_m=\pi_{n+1}^*\kappa_m+\psi_{n+1}^m
\end{equation}
and the straightforward relation
\begin{equation}  \label{eq:pbpsi}
\psi_j=\pi_{n+1}^*\psi_j,\ j=1,...,n.
\end{equation}

To any $\Sigma$ bundle $\pi:X\to B$ with $n$ sections $s_i$ corresponds the pull-back $\Sigma$ bundle $\pi^*X\to X$ with $n$ sections $\pi^*s_i$ and a further tautological section $s_{n+1}$.
In some sense the section $s_{n+1}$ gives all possible ways to add an $(n+1)$st section to the bundle over $B$.  In this context the forgetful map has two interpretations.  As the map $\pi_{n+1}^*:H^*(B)\to H^*(B)$ discussed above, and also as $\pi^*:H^*(B)\to H^*(X)$.  The two are related by 
\[ s_{n+1}^*\circ\pi^*=\pi_{n+1}^*.\]
The pull-back relation (\ref{eq:pbkap}) looks the same for $\pi^*$
\renewcommand{\theequation}{\arabic{equation}a}
\addtocounter{equation}{-2}
\begin{equation}  \label{eq:pbkapa}
\kappa_m=\pi^*\kappa_m+\psi_{n+1}^m
\end{equation}
whereas the relation (\ref{eq:pbpsi}) needs to be adjusted  to
\begin{equation} \label{eq:pbpsia}
\psi_j=\pi^*\psi_j+PD[s_j(B)],\quad j=1,...,n.
\end{equation}
\renewcommand{\theequation}{\arabic{equation}}

We have yet to mention that the tautological section $s_{n+1}:X\to\pi^*X$ does not have disjoint image from the other sections $\pi^*s_i$.  After blowing up to separate the images of the sections, we are naturally led to consider surface bundles $\pi:X\to B$ that allow fibres with mild singularities.  More precisely, the singular fibres may be {\em stable curves} - they consist of  a collection of smooth components meeting at nodal singularities with the property that each component has multiplicity 1, and negative Euler characteristic after we subtract all labeled points and common points with other components.   We call $X$ a bundle of stable curves or simply a bundle with singular fibres, although strictly it is no longer a fibre bundle.  The cohomology classes $\psi_i$ and $\kappa_m$ extend to this situation.  Their definitions are best understood when we put a continuous family of conformal structures on the fibres, or we assume the stronger property that $X$ and $B$ are complex analytic varieties.  Define $\gamma=K_X\otimes\pi^*K_B^{-1}$,  essentially the vertical canonical bundle (relative dualising sheaf.)  This coincides with the definition above on smooth fibres and generalises the definition to singular fibres.  One can make sense of sections of this bundle along singular fibres in terms of meromorphic 1-forms with simple poles and conditions on residues \cite{HMoMod} but we will not explain this here.  The definitions of $\psi_i$ and $\kappa_m$ are as above.  Relations (\ref{eq:pbkapa}) and (\ref{eq:pbpsia}) generalise to bundles of stable curves.  Proofs can be found in \cite{ACoCom} and  \cite{WitTwo}.

A simple example will demonstrate the definitions and relations.  Let $X$ be the blow-up of $\bp^1\times\bp^1$ at the three points $(0,0)$, $(1,1)$ and $(\infty,\infty)$.  The map from $X$ to the first $\bp^1$ factor realises $X$ as a surface bundle
\[ \begin{array}{cccll}\bp^1&\hookrightarrow&\ \ X&\\&&\pi\downarrow&\uparrow s_i\\&&\quad \bp^1&\end{array}
\quad\quad i=1,...,4\]
which we equip with four sections $s_1(z)=(z,0)$, $s_2(z)=(z,1)$, $s_3(z)=(z,\infty)$ and $s_4(z)=(z,z)$.  The general fibre is genus 0 with 4 labeled points, and the singular fibres, at $0$, $1$ and $\infty$, are stable curves with two irreducible components each with two labeled points (and a common intersection point.)  We can generate $H_2(X)$ by $H$, $F$, $E_1$, $E_2$ and $E_3$ where $E_i$ are the exceptional divisors of the blow-up and $H=\bp^1\times\{w\}$ and $F=\{z\}\times\bp^1$ for any $w$ and $z$ different from 0, 1 and $\infty$.  We use these curves to represent their divisor class, homology class and their Poincare dual cohomology class.  Then 
\[ c_1(\gamma)=-2H+E_1+E_2+E_3\quad\Rightarrow \tilde{\kappa}_1=c_1(\gamma)^2=-3,\]
\[c_1\left(\gamma\left[\sum s_i(B)\right]\right)=2H+F-E_1-E_2-E_3\quad\Rightarrow\kappa_1=
c_1\left(\gamma\left[\sum s_i(B)\right]\right)^2=1,\]
\[\psi_1=c_1(\gamma)\cdot s_1(B)=c_1(\gamma)\cdot(H-E_1)=1=\psi_i,\quad i=2,3,4.\]
Since $X$ is the blow-up of the pull-back of the $\bp^1$ bundle over a point with three sections, (\ref{eq:pbkapa}) and (\ref{eq:pbpsia}) are also evident.

\subsection{Intersection numbers}

Let us use $\modm_{g,n}$ to notate the moduli space of genus $g$ curves with $n$ labeled points, which is isomorphic to the moduli space of genus $g$ hyperbolic surfaces with $n$ labeled cusps, $\modm_{g,n}({\bf L})$ with ${\bf L}={\bf 0}$, and $\overline{\modm}_{g,n}$ the Deligne-Mumford compactification which adds stable curves to $\modm_{g,n}$.  Wolpert \cite{WolWei} showed that the symplectic structure $\omega$ on $\modm_{g,n}$ extends to $\overline{\modm}_{g,n}$.  The $\psi_i$ and $\kappa_m$ classes naturally live in $H^*(\overline{\modm}_{g,n})$.  They are associated to a universal surface bundle over $\overline{\modm}_{g,n}$, essentially given by $\overline{\modm}_{g,n+1}$ with map forgetting the last labeled point, and any bundle $X$ equipped with conformal structures on fibres is the pull-back of the universal bundle under a map $B\to\overline{\modm}_{g,n}$.

\begin{theorem}[Mirzakhani]  \label{th:coeff}
The coefficient $C_{\alpha}$ of $L_1^{2\alpha_1}...L_n^{2\alpha_n}$ in $V_{g,n}({\bf L})$ is
\begin{equation}  \label{eq:coeff}
C_{\alpha}=\frac{1}{2^{|\alpha|}\alpha!(3g-3+n-|\alpha|)!}\int_{\overline{\modm}_{g,n}}\psi_1^{\alpha_1}...\psi_n^{\alpha_n}\omega^{3g-3+n-|\alpha|}
\end{equation}
\end{theorem}
This is proven in \cite{MirWei} by showing that $\modm_{g,n}({\bf L})$ is the symplectic quotient of a larger symplectic manifold by a Hamiltonian $T^n$ action, where a fixed value of the moment map corresponds to fixing the lengths $L_1,...,L_n$ of the geodesic boundary components.  Any such quotient is equipped with $n$ line bundles coming from the $T^n$ action, and their Chern classes are related to the coefficients of the volume polynomial.  In \cite{MirWei} Mirzakhani used this together with her recursion relation for the volume polynomials to give a new proof of Witten's conjecture \cite{WitTwo} regarding intersections of $\psi$ classes on $\overline{\modm}_{g,n}$.  In the original proof of Witten's conjecture, Kontsevich \cite{KonInt} calculated the Laplace transform of the top degree terms of $V_{g,n}({\bf L})$.  It would be interesting to understand the Laplace transform of the whole polynomial $V_{g,n}({\bf L})$.

In the following, write $\psi^{\alpha}$ for $\psi_1^{\alpha_1}...\psi_n^{\alpha_n}$ and ignore the term if there is an $\alpha_j<0$.  For ease of reading, note that in all formulae the variable $j$ sums from $0$ to $m$ while the variable $k$ sums from $1$ to $n$.
\begin{lemma}
The equation
\[ V_{g,n+1}({\bf L},2\pi i)=\sum_{k=1}^n\int_0^{L_k}L_kV_{g,n}({\bf L})dL_k\]
is equivalent to
\begin{equation}   \label{eq:string2}
\sum_{j=0}^m(-1)^j{m\choose j}\int_{\overline{\modm}_{g,n+1}}\psi^{\alpha}\psi_{n+1}^j\kappa_1^{m-j}= \sum_{k=1}^n\int_{\overline{\modm}_{g,n}}\psi_1^{\alpha_1}..\psi_k^{\alpha_k-1}..\psi_n^{\alpha_n}
 \kappa_1^m
 \end{equation}
 for all $\alpha$ and $m$.
\end{lemma}
\begin{proof}
Assume that $|\alpha|+m=3g-2+n$ since otherwise (\ref{eq:string2}) is zero on both sides.
By (\ref{eq:coeff}) and substitution of $L_{n+1}^{2j}$ with $(2\pi i)^{2j}$, the coefficient of $L_1^{2\alpha_1}...L_n^{2\alpha_n}$ in $V_{g,n+1}({\bf L},2\pi i)$ is
\[\sum_{j=0}^m\frac{(2\pi i)^{2j}}{2^{|\alpha|+j}\alpha!j!(m-j)!}\int_{\overline{\modm}_{g,n+1}}\psi^{\alpha}\psi_{n+1}^j\omega^{m-j}\quad\quad\quad\quad\quad\quad \]
\[ =\sum_{j=0}^m\frac{(2\pi i)^{2j}}{2^{|\alpha|+j}\alpha!j!(m-j)!}\int_{\overline{\modm}_{g,n+1}}\psi^{\alpha}\psi_{n+1}^j(2\pi^2\kappa_1)^{m-j}\]
\[ =\frac{2^{m-|\alpha|}\pi^{2m}}{\alpha!\ m!}\sum_{j=0}^m(-1)^j{m\choose j}\int_{\overline{\modm}_{g,n+1}}\psi^{\alpha}\psi_{n+1}^j\kappa_1^{m-j}\]
where we have used the identity $\omega=2\pi^2\kappa_1$ proven in \cite{WolHom}.

The coefficient of $L_1^{2\alpha_1}...L_n^{2\alpha_n}$ in $\int_0^{L_k}L_kV_{g,n}dL_k$ is
\[\frac{\alpha_k}{2^{|\alpha|-1}(2\alpha_k)\alpha!m!}\int_{\overline{\modm}_{g,n}}\psi_1^{\alpha_1}..\psi_k^{\alpha_k-1}..\psi_n^{\alpha_n}\omega^m\quad\quad\quad\quad\quad\quad \]
\[ =\frac{2^{m-|\alpha|}\pi^{2m}}{\alpha!\ m!}\int_{\overline{\modm}_{g,n}}\psi_1^{\alpha_1}..\psi_k^{\alpha_k-1}..\psi_n^{\alpha_n}\kappa_1^m.\]
Add this expression over $k=1,...,n$ and divide both sides by the factor\\ $2^{m-|\alpha|}\pi^{2m}/\alpha!m!$ to prove the lemma.
\end{proof}

\begin{lemma}
The equation
\[ \frac{\partial V_{g,n+1}}{\partial L_{n+1}}({\bf L},2\pi i)=2\pi i(2g-2+n)V_{g,n}({\bf L})\]
is equivalent to
\begin{equation} \label{eq:dilaton2}
 \sum_{j=0}^m(-1)^j{m\choose j}\int_{\overline{\modm}_{g,n+1}}\psi^{\alpha}\psi_{n+1}^{j+1}\kappa_1^{m-j}= (2g-2+n)\int_{\overline{\modm}_{g,n}}\psi^{\alpha}\kappa_1^m
 \end{equation}
\end{lemma}
\begin{proof}
The proof is much like the proof of the previous lemma.  The coefficient of $L_1^{2\alpha_1}...L_n^{2\alpha_n}$ in $\partial V_{g,n+1}/\partial L_{n+1}({\bf L},2\pi i)$ is
\[\sum_{j=0}^m\frac{(2j+2)(2\pi i)^{2j+1}}{2^{|\alpha|+j+1}\alpha!(j+1)!(m-j)!}\int_{\overline{\modm}_{g,n+1}}\psi^{\alpha}\psi_{n+1}^{j+1}\omega^{m-j}\quad\quad\quad\quad\quad\quad  \]
\[ =2\pi i\frac{2^{m-|\alpha|}\pi^{2m}}{\alpha!\ m!}\sum_{j=0}^m(-1)^j{m\choose j}\int_{\overline{\modm}_{g,n+1}}\psi^{\alpha}\psi_{n+1}^{j+1}\kappa_1^{m-j}\]
and the coefficient of $L_1^{2\alpha_1}...L_n^{2\alpha_n}$ in $V_{g,n}$ is
\[\frac{2^{m-|\alpha|}\pi^{2m}}{\alpha!\ m!}\int_{\overline{\modm}_{g,n}}\psi^{\alpha}\kappa_1^m\]
so the equivalence follows.

\end{proof}

\begin{proof} [Completion of the proof of Theorem~\ref{th:string}]
It suffices to prove the relations (\ref{eq:string2}) and (\ref{eq:dilaton2}).  Notice that when $m=0$, 
(\ref{eq:string2}) and (\ref{eq:dilaton2}) are respectively the string and dilaton equations which were proven by Witten in \cite{WitTwo}.  The method of proof for the more general identities is similar.

Let $\pi:X\to B$ be a bundle of  stable curves with $n$ disjoint sections and $\pi^*X$ the pull-back bundle with $n+1$ sections.  Blow up $\pi^*X$ along the intersections of images of sections to get a bundle of  stable curves over $X$ with $n+1$ disjoint sections $s_1,...,s_{n+1}$.  Our aim is to compare $\psi$ and $\kappa$ classes in $H^*(X)$ and $H^*(B)$.

Take the integrand of the left hand side of (\ref{eq:string2}) and consider its image under the umkehr map.
\[
\pi_!\left\{\sum_{j=0}^m(-1)^j{m\choose j}\psi_{n+1}^j\kappa_1^{m-j}\prod_{k=1}^n\psi_k^{\alpha_k}\right\}= 
\pi_!\left\{(\kappa_1-\psi_{n+1})^m\prod_{k=1}^n\psi_k^{\alpha_k}\right\}\quad\quad\quad\quad\quad\]
\[\quad\quad\quad\quad\quad\quad\quad\quad\quad\quad\quad=\pi_!\left\{(\pi^*\kappa_1^m)\prod_{k=1}^n\left(\pi^*\psi_k^{\alpha_k}+PD[s_k(B)]\cdot\pi^*\psi_k^{\alpha_k-1}\right)\right\}\]
\[=\kappa_1^m\sum_{k=1}^n\psi_1^{\alpha_1}..\psi_k^{\alpha_k-1}..\psi_n^{\alpha_n}. \]
To get from the first line to the second line we have used the pull-back formulae  (\ref{eq:pbkapa}) and  (\ref{eq:pbpsia}) and the fact that $\pi^*$ is a ring homomorphism, so in particular $(\pi^*\eta)^m=\pi^*(\eta^m)$.   
To get from the second line to the third line we have used the fact that $\pi_!:H^*(X)\to H^*(B)$ is an $H^*(B)$ module homomorphism, i.e. $\pi_!(\xi\pi^*\eta)=\pi_!(\xi)\eta$, together with the explicit evaluations
\[\pi_!(1)=0,\quad\pi_!(s_i(B))=1\]
most easily calculated from the Poincare duality description of $\pi_!$.  Thus, in the product the image under $\pi_!$ of the highest degree term $\pi^*(\kappa_1^m\psi^{\alpha})$ is zero, the image of the second highest degree term constitutes the expression in the third line, and the lower order terms vanish since they contain products $PD[s_j(B)]\cdot PD[s_k(B)]=0$ because the images of $s_j$ and $s_k$ are disjoint.

Since
\[\int_X\eta= \int_B\pi_!\eta\]
choose $X=\overline{\modm}_{g,n+1}$ and $B=\overline{\modm}_{g,n}$, so (\ref{eq:string2}) follows.\\ 

The proof of (\ref{eq:dilaton2}) is similar.  Again apply the umkehr map to the integrand of the left hand side of (\ref{eq:dilaton2}).
\[
\pi_!\left\{\sum_{j=0}^m(-1)^j{m\choose j}\psi_{n+1}^{j+1}\kappa_1^{m-j}\prod_{k=1}^n\psi_k^{\alpha_k}\right\}= 
\pi_!\left\{\psi_{n+1}\cdot(\kappa_1-\psi_{n+1})^m\prod_{k=1}^n\psi_k^{\alpha_k}\right\}\quad\quad\]
\[\quad\quad\quad\quad\quad\quad\quad\quad\quad\quad=\pi_!\left\{\psi_{n+1}\cdot(\pi^*\kappa_1^m)\prod_{k=1}^n\left(\pi^*\psi_k^{\alpha_k}+PD[s_k(B)]\cdot\pi^*\psi_k^{\alpha_k-1}\right)\right\}\]
\[=(2g-2+n)\kappa_1^m\prod_{k=1}^n\psi_k^{\alpha_k}. \]
To go from the second line to the third line note that $\psi_{n+1}$ coincides with the twisted Euler class $e_{n+1}$ that satisfies (\ref{eq:euler}) and hence $\pi_!\psi_{n+1}=2g-2+n$ and $\psi_{n+1}\cdot PD[s_k(B)]=0$.  Thus the top degree term constitutes the expression in the third line, and all lower degree terms vanish.
\end{proof}

Equations (\ref{eq:string}) and (\ref{eq:dilaton}) suggest that a direct analysis of the moduli space of cone surfaces with cone angle $\theta\approx 2\pi$, or more accurately an infinitesimal analysis near $2\pi$, will gives rise to intriguing phenomena.  Equation (\ref{eq:dilaton}) seems plausible since the removed cone point is free to wander around each hyperbolic surface with area $2\pi(2g-2+n)$, so the change in volume is related to integrating over the smaller moduli space and along each fibre.  Intuition for equation (\ref{eq:string}) seems less obvious.

\section{Use of recursion relations}  \label{sec:rec}
\subsection{Classical volumes and low genus algorithms.}
The relations (\ref{eq:string}) and (\ref{eq:dilaton}) can be used to give convenient information about the volumes of moduli spaces.  In the following, when there is only one boundary component we use the variable $L=L_1$.
\begin{enumerate}
\item[(i)] When there is exactly one marked point  the volume factorises
\begin{equation}  \label{eq:factor}
V_{g,1}(L)=(L^2+4\pi^2)P_g(L)
\end{equation}
\item[(ii)]  The classical volumes of moduli spaces, i.e. when there are no marked points, are encoded in Mirzakhani's volume polynomials.  In terms of the polynomial $P_g$ defined by (\ref{eq:factor})
\begin{equation} \label{eq:vol0}
\left.V_{g,0}=\frac{P_g(2\pi i)}{g-1}=\frac{1}{(g-1)}\frac{V_{g,1}}{L^2+4\pi^2}\right|_{L=2\pi i}
\end{equation}
\item[(iii)]  The relation (\ref{eq:string}) uniquely determines $V_{0,n+1}$ from $V_{0,n}$.
\item[(iv)]  The relations (\ref{eq:string}) and (\ref{eq:dilaton}) uniquely determine $V_{1,n+1}$ from $V_{1,n}$.
\end{enumerate}
Equations (\ref{eq:factor}) and (\ref{eq:vol0}) follow from the proof of Theorem~\ref{th:string}.  In that proof we see that when $n=0$ the right hand side of (\ref{eq:string}) must be zero, which can be interpreted as there being no boundary lengths to integrate, and hence $V_{g,1}(L)$ possesses a factor of $L^2+4\pi^2$.  To get (\ref{eq:vol0}), use (\ref{eq:dilaton}) together with $V_{g,1}(2\pi i)=0$:
\begin{eqnarray*}
 2\pi i(2g-2)V_{g,0}=\left. \frac{dV_{g,1}}{dL}\right|_{L=2\pi i}
&=& \lim_{L\to 2\pi i}\frac{V_{g,1}(L)}{L-2\pi i}\\&=& \lim_{L\to 2\pi i}\frac{4\pi iV_{g,1}(L)}{L^2+4\pi^2}=4\pi iP_g(2\pi i)
\end{eqnarray*}
\begin{proof}[Proof of (iii).]
This  follows from the elementary fact that a symmetric polynomial $f(x_1,...,x_n)$ of degree less than $n$ is uniquely determined by evaluation of one variable at any $a\in\bc$, $f(x_1,...,x_{n-1},a)$.  To see this, suppose otherwise.  Any symmetric $g(x_1,...,x_n)$ of degree less than $n$ that evaluates at $a$ as $f$ does, satisfies 
\begin{eqnarray*}
f(x_1,...,x_{n-1},a)-g(x_1,...,x_{n-1},a)&=&(x_n-a)P(x_1,...,x_n)\\
&=&Q(x_1,...,x_n)\prod_{j=1}^n(x_j-a)
\end{eqnarray*}
but the degree is less than $n$ so the difference is identically 0.

The volume $V_{0,n+1}$ is a symmetric degree $n-2$ polynomial in $L_1^2,...,L_{n+1}^2$ so it is uniquely determined by evaluation at $L_{n+1}=2\pi i$, and this is determined by $V_{0,n}$ via (\ref{eq:string}).
\end{proof}
The proof of (iv) is similar to the proof of (iii).  The degree of $V_{1,n+1}$ is equal to $n+1$ so the proof of (iii) shows that (2) determines $V_{1,n+1}$ from $V_{1,n}$ up to the constant $c$ in $V_{1,n+1}+c\prod_{j=1}^n(L_j^2+4\pi^2)$.  Now use (\ref{eq:dilaton}) to determine $c$, and hence $V_{1,n+1}$.\\

Statement (iii) can be converted to an algorithm for calculating $V_{0,n}$.  The algorithm using (\ref{eq:string}) turns out to be much more efficient than the algorithm coming from Mirzakhani's relation (\ref{eq:rec}) in genus 0, which needs $V_{0,n-1}$ and pairs $V_{0,n_1}, V_{0,n_2}$ for all $n_1+n_2=n+1$, to produce $V_{0,n}$.   We have included a simple MAPLE routine in the appendix for calculating $V_{0,n}$ using (\ref{eq:string}).  (The notion of a ``more efficient" algorithm is not so precise here.  We have merely compared the speeds of different calculations on MAPLE.)  

In genus 0, the string equation - (\ref{eq:string2}) with $m=0$ - leads to an explicit formula for the top coefficients, or equivalently the following formula for genus 0 intersection numbers without kappa classes:
\[\int_{\overline{\modm}_{g,n}}\psi_1^{\alpha_1}...\psi_n^{\alpha_n}={n-3\choose\alpha_1,...,\alpha_n}.\]
It seems reasonable to guess that when $g=0$ the relation (\ref{eq:string2}) leads to an explicit combinatorial description of all genus 0 intersection numbers with powers of $\kappa_1$, or equivalently all coefficients of $V_{0,n}$.  Zograf \cite{ZogWei} has recursion relations between the constant coefficients $V_{0,n}({\bf 0})$.

\subsection{Higher derivatives}
We expect to have expressions for higher derivatives $\partial^kV_{g,n+1}/\partial L_{n+1}^k$ evaluated at $L_{n+1}=2\pi i$.  Evidence comes from the fact that (\ref{eq:string}) and (\ref{eq:dilaton}) use generalised versions of the string and dilaton equations.  The Virasoro relations are a sequence of relations for the top degree terms of $V_{g,n}({\bf L})$, with first two relations in the sequence the string and dilaton equations, so may also have versions in terms of evaluations of derivatives  of the volume polynomial at $L_{n+1}=2\pi i$.  The Virasoro relations recursively determine the top degree coefficients of the volume polynomials by using the relations in a clever way.   In recent work \cite{MSaMir}, Mulase  and Safnuk showed how to extend the Virasoro relations to the full volume polynomials.  It would be desirable to instead determine the polynomials recursively by relying on the more familiar fact that the derivatives of a function evaluated at a point determine the function. It would be interesting to know if one can express the results \cite{MSaMir} in terms of derivatives  of the volume polynomial at $L_{n+1}=2\pi i$.

In principle, we can use Mirzakhani's recursion relation to get expressions for higher derivatives of the volume evaluated at $L_{n+1}=2\pi i$.  Differentiate the equation
\[
\frac{\partial (L_{n+1}V_{g,n+1})}{\partial L_{n+1}}=\ay_{g,n+1}+\bb_{g,n+1}.
\]
to get
\[\frac{\partial^2( L_{n+1}V_{g,n+1})}{\partial L_{n+1}^2}=\frac{\partial \ay_{g,n+1}}{\partial L_{n+1}}+\frac{\partial \bb_{g,n+1}}{\partial L_{n+1}}\]
and evaluate at $L_{n+1}$.   Substitute the equation for the first derivative, to get 
the following equation for the second derivative.  Put $\mathcal{E}=\sum_{j=1}^nL_j\partial/\partial L_j$, the Euler vector field:
\[
\frac{\partial^2 V_{g,n+1}}{\partial L_{n+1}^2}({\bf L},2\pi i)=\mathcal{E}\cdot V_{g,n}({\bf L})-2\pi i(4g-4+n)V_{g,n}({\bf L}).
\]
By taking higher derivatives of Mirzakhani's relation we can recursively get equations for higher derivatives.  The strength of (\ref{eq:string}) and (\ref{eq:dilaton}) is the simplification of Mirzakhani's relations (\ref{eq:rec}).  It is not clear that the higher derivative relations obtained by the method above possess this same strength.

\section*{Appendix}

{\sc MAPLE routine for calculating $V_{0,n}$.}

\begin{mapleinput}
\mapleinline{active}{1d}{\begin{Maple Normal}\QTR{Maple Maple Input}{}\QTR{Maple Maple Input}
{#  input: symmetric polynomial f in n variables L1,...,Ln} \QTR
{Maple Maple Input}{
} \QTR{Maple Maple Input} {
\newline # output: symmetric polynomial S in n+1 variables L1,...,L(n+1)  
\newline # satisfying S(L(n+1)=0)=f}\QTR{Maple Maple Input}{
}\QTR{Maple Maple Input}{
\newline sym:=proc(f) local i,j,k,m,S,T,T1,prod,sum,epsilon:  
\newline S:=f: }\QTR{Maple Maple Input}{
}\QTR{Maple Maple Input}{
\newline epsilon:=array[1..100]: 
\newline for i from 1 to 100 do epsilon[i]:=0: od: }\QTR{Maple Maple Input}{
}\QTR{Maple Maple Input}{
\newline while epsilon[n+1]<1 do}\QTR{Maple Maple Input}{
}\QTR{Maple Maple Input}{
\newline T:=subs({seq(L||j=(1-epsilon[j])*L||j,j=1..n)},f): 
\newline T1:=0:}\QTR{Maple Maple Input}{
 }\QTR{Maple Maple Input}{
\newline for i from 1 to n do}\QTR{Maple Maple Input}{
\newline prod:=1:
 \newline  for j from i+1 to n+1 do 
 \newline prod:=prod*(1-epsilon[j]) 
 \newline od:}\QTR{Maple Maple Input}{
 \newline T1:=T1+prod*subs(L||i=L||(n+1),T): }\QTR{Maple Maple Input}{
\newline od:}\QTR{Maple Maple Input}{
}\QTR{Maple Maple Input}{
\newline sum:=0: for k from 1 to n do sum:=sum+epsilon[k] od: }\QTR{Maple Maple Input}{
}\QTR{Maple Maple Input}{
\newline S:=S+(-1)\symbol{94}sum*T1: }\QTR{Maple Maple Input}{
}\QTR{Maple Maple Input}{
\newline for k from 1 to 100 do}\QTR{Maple Maple Input}{
\newline  if epsilon[k]=1 then epsilon[k]:=0 
\newline else epsilon[k]:=1: k:=100 end if:}\QTR{Maple Maple Input}{
\newline od:}\QTR{Maple Maple Input}{
 }\QTR{Maple Maple Input}{
\newline od:}\QTR{Maple Maple Input}{
\newline S:=simplify(S):}\QTR{Maple Maple Input}{
\newline end:}\end{Maple Normal}}{}
\end{mapleinput}
\begin{mapleinput}
\mapleinline{active}{1d}{\begin{Maple Normal}\QTR{Maple Maple Input}{}\QTR{Maple Maple Input}{# calculate the genus zero volumes recursively from evaluation 
\newline # of V_{(0,n+1)} at L(n+1)=2*Pi*I}\QTR{Maple Maple Input}{
\newline for n from 3 to 12 do}\QTR{Maple Maple Input}{
 \newline P:=0: 
 \newline for j from 1 to n do 
 \newline P:=P+int(L||j*V[n],L||j) 
 \newline od:}\QTR{Maple Maple Input}{
 \newline Q0:=P: 
 \newline C0:=simplify(coeff(Q0,Pi,0)): 
 \newline sim:=sym(C0):
 \newline V[n+1]:=sim:}\QTR{Maple Maple Input}{
\newline for k from 1 to n-2 do}\QTR{Maple Maple Input}{
\newline  P||k:=sim-C||(k-1):}\QTR{Maple Maple Input}{
 \newline Q||k:=subs(L||(n+1)=2*Pi*I,Q||(k-1)-P||k*Pi\symbol{94}(2*k-2)):}\QTR{Maple Maple Input}{
\newline  C||k:=simplify(coeff(Q||k,Pi,2*k)):}\QTR{Maple Maple Input}{
\newline  sim:=sym(C||k):}\QTR{Maple Maple Input}{
\newline  V[n+1]:=V[n+1]+sim*Pi\symbol{94}(2*k):}\QTR{Maple Maple Input}{
\newline od:}\QTR{Maple Maple Input}{
\newline od:}\end{Maple Normal}}{}
\end{mapleinput}


\begin{thebibliography}{99}

\bibitem{ACoCom} Arbarello, Enrico and Cornalba, Maurizio 
\emph{Combinatorial and algebro-geometric cohomology classes on the moduli spaces of curves.}
 J. Algebraic Geom. {\bf 5} (1996), 705--749.

\bibitem{HMoMod} Harris, Joe and Morrison, Ian 
\emph{Moduli of curves.} Graduate Texts in Mathematics {\bf 187}. Springer-Verlag, New York, 1998. 

\bibitem{KonInt}  Kontsevich, Maxim 
\emph{Intersection theory on the moduli space of curves and the matrix Airy function.} 
Communications in Mathematical Physics, {\bf 147} (1992), 1Ð23.

\bibitem{KMroEmb}
Kronheimer, P. B. and Mrowka, T. S.
\emph{Embedded surfaces and the structure of Donaldson's polynomial invariants. }
J. Diff. Geom. {\bf 41} (1995), 573--734.

\bibitem{McSRem} 
McShane, Greg
\emph{A remarkable identity for lengths of curves.} Ph.D. Thesis, University of
Warwick, 1991.

\bibitem{MirSim}
Mirzakhani, Maryam 
\emph{Simple geodesics and Weil-Petersson volumes of moduli spaces of bordered Riemann surfaces.}Preprint.

\bibitem{MirWei}
Mirzakhani, Maryam 
\emph{Weil-Petersson volumes and intersection theory on the moduli space of curves.}
Preprint.

\bibitem{MorCha} Morita, Shigeyuki 
\emph{Characteristic classes of surface bundles.} 
Invent. Math. {\bf 90} (1987), 551--577. 

\bibitem{MSaMir} Mulase, Motohico  and Safnuk, Brad 
\emph{Mirzakhani's recursion relations, Virasoro constraints and the KdV hierarchy.}
math.QA/0601194 

\bibitem{TWZGen} Ser Peow Tan, Yan Loi Wong, Ying Zhang
\emph{Generalizations of McShane's identity to hyperbolic cone-surfaces.}
To appear, J. Diff. Geom. {\bf math.GT/0404226.}

\bibitem{WitTwo} Witten, Edward
\emph{Two-dimensional gravity and intersection theory on moduli space.} Surveys in differential geometry (Cambridge, MA, 1990), 243--310, Lehigh Univ., Bethlehem, PA, 1991.

\bibitem{WolHom} Wolpert, Scott 
\emph{On the homology of the moduli space of stable curves. }
Ann. of Math. (2) {\bf 118} (1983), 491--523.

\bibitem{WolWei} Wolpert, Scott 
\emph{On the Weil-Petersson geometry of the moduli space of curves.}
Amer. J. Math. {\bf 107} (1985), 969--997.

\bibitem{ZogWei}  Zograf, Peter 
\emph{The Weil-Petersson volume of the moduli space of punctured spheres.} 
Mapping class groups and moduli spaces of Riemann surfaces, 367--372, Contemp. Math., 150, Amer. Math. Soc., Providence, RI, 1993.

\end{thebibliography}
\end{document}